\documentclass[11pt]{article}
\usepackage{graphicx,amsthm}

\newcommand{\abs}[1]{\left|#1\right|}
\newtheorem{theorem}{Theorem}

\begin{document}

\title{The Imbedding Sum of a Graph}
\author{Robert G. Rieper\\William Paterson University}
\maketitle
\abstract{The automorphisms of a graph act naturally on its set of labeled imbeddings to produce its
unlabeled imbeddings. The imbedding sum of a graph is a polynomial that contains useful information about a graph's
labeled and unlabeled imbeddings.  In particular, the polynomial enumerates the number of different ways the unlabeled
imbeddings can be vertex colored and enumerates the labeled and unlabeled imbeddings by their symmetries.}

\section{Introduction}
\subsection{Motivation}
Graphical enumeration is a well-established discipline that counts the number of graphs with a selected list of
properties. These counting problems frequently arise in the structure of chemical compounds, in computer science and in
combinatorial designs. In the article by Mull, Rieper, and White~\cite{MRW88} a technique was developed for counting
the number of different imbeddings a graph has on orientable surfaces. The result was extended to include imbeddings on
nonorientable surfaces by Kwak and Lee~\cite{KL94}. For example, the complete graph on four vertices has three
different imbeddings among the orientable surfaces (one on the plane and two on the torus) and eleven among all
surfaces including the nonorientable cases.  Similar results for complete bipartite graphs were found by
Mull~\cite{M99}.

The purpose of this article is to show how these results can be extended to obtain information about graph imbeddings
beyond just their quantity. We show, for example, that the three different orientable imbeddings of $K_4$ become five
if two of the vertices are distinguished from the other two by coloring them black and the others white.  We also show
how the imbeddings are distributed according to the symmetries that they have.

\subsection{Graph imbeddings}
An \emph{imbedding} of a graph $G$, considered to be a simplicial complex, into a closed orientable 2-manifold $S$, is
a homeomorphism $i$ of $G$ into $S$, $i:G\rightarrow S$.  If every component of $S-i(G)$ is a 2-cell (a homeomorph of
an open disk), then $i$ is a \emph{2-cell } imbedding.  If an orientation is provided for $S$, then $i$ is an
\emph{oriented} imbedding.  The orientable 2-manifolds considered here are the connected compact topological spaces
each of which is characterized as a sphere with handles.  Such a space is called an \emph{orientable surface} and the
number of handles is its \emph{genus}.  We regard two oriented imbeddings $i:G \rightarrow S$ and $j:G \rightarrow S$
as the same if there exists an orientation-preserving homeomorphism $h$ of $S$ onto $S$ such that $h \circ i=j$.
Henceforth, a 2-cell oriented imbedding is referred to as a \emph{labeled imbedding}.

Each labeled imbedding of a graph can be described by the following well-known scheme (see Gross and
Tucker~\cite{GT87}). In a small neighborhood of a vertex $v$ we observe the counterclockwise cyclic order of the edges
incident with $v$. If the graph has no loops or multiple edges, then we record the vertices adjacent to $v$ in this
order as the \emph{rotation} at $v$. The \emph{rotation system} is the vertex-indexed set of these rotations. If the
vertex $v$ is adjacent to $\mathrm{degree}(v)$ other vertices, then there are $(\mathrm{degree}(v)-1)!$ different
rotations at $v$. The product over the vertex set of these numbers is the quantity of different labeled imbeddings.

For a given graph, many of its labeled imbeddings resemble one another.
 Formally, two labeled imbeddings $i:G \rightarrow S$
and $j:G \rightarrow S$ are \emph{congruent} if there exists a graph automorphism $h:G \rightarrow G$ such that $i
\circ h=j$. A congruence class is called an \emph{unlabeled imbedding}. Informally, a labeled imbedding of a graph is a
drawing of the graph on a surface where each vertex receives a label. In an unlabeled imbedding the vertex labels are
omitted. We remark that the definitions of rotation scheme and of labeled and unlabeled imbeddings can be naturally
extended to include graphs with multiple edges or loops and to directed graphs.


\section{The Imbedding Sum of a Graph}
\label{section2} Let $G$ be a connected \emph{simple} graph with vertex set $V$ where a simple graph has no loops or
multiple edges. The restriction to simple graphs enables us to more easily describe the automorphisms (they act on the
vertices alone) and is easily circumvented by modeling a more general graph with a suitable simple one. The
\emph{neighborhood} $N(v)$ of the vertex $v$ is the set of vertices adjacent to $v$. Each rotation at $v$ is a cyclic
permutation of the neighborhood of $v$, denoted $\rho_v:N(v) \rightarrow N(v)$. The rotation system is the indexed set
$\rho = \{\rho_v\}_{v \in V}$. A \emph{map} is a pair $M=(G,\rho)$ and provides an algebraic correspondence with the
labeled imbeddings.

We denote the automorphism group of $G$ by $\Gamma(G)$ (or simply $\Gamma$ if $G$ is understood) and define two maps
$M=(G,\rho)$ and $M'=(G,\rho')$ to be \emph{equivalent} if there exists an automorphism $\gamma \in \Gamma (G)$ such
that $\gamma \rho_v \gamma^{-1}=\rho_{\gamma (v)}^{'}$ for all vertices $v$.  The graph automorphism $\gamma$ can be
interpreted as relabeling the vertices of $G$. Equivalent maps correspond with congruent labeled imbeddings. A
\emph{map automorphism} for $M=(G,\rho)$ is a graph automorphism giving $M$ equivalent to itself. The set of map
automorphisms for $M=(G,\rho )$ form the map-automorphism subgroup $\Gamma_M(G)$ of $\Gamma(G)$. Thus, the
graph-automorphism group $\Gamma(G)$ acts on the set of maps and $\Gamma_M(G)$ is the stabilizer of the map $M$ under
this action. We therefore have the following result (Biggs \cite{B71}).

\begin{theorem}
\label{thm:index} The number of labeled imbeddings of the graph $G$ in the unlabeled class containing the map $M$ is
the index of $\Gamma_M(G)$ in $\Gamma(G)$.
\end{theorem}

The next three theorems, called the \emph{counting theorems} are those derived in \cite{MRW88} and are to be used
extensively in the applications to follow. For this reason they are reproduced here without proof. For each graph
automorphism $\gamma \in \Gamma$ we let $F(\gamma)$ be the set of maps which are fixed by $\gamma$. That is, $M \in
F(\gamma)$ if and only if $M$ is $\gamma$-equivalent to itself. An application of Burnside's lemma yields the next
result.

\begin{theorem}
\label{thm:Burnside} The number of different unlabeled imbeddings of a graph is $\frac{1}{\abs{\Gamma}}\sum_{\gamma \in
\Gamma} \abs{F(\gamma)}$.
\end{theorem}

The cardinality of a fixed set is determined as follows. Each automorphism $\gamma$ of the graph $G$ does double duty.
It acts on the set of vertices and the set of maps. If $v$ is a vertex of $G$, then we let $l(v,\gamma)$ denote the
length of the orbit which contains the vertex $v$ under the action of $\gamma$. We define the \emph{fixed set at $v$},
denoted $F_v(\gamma)$, to be the set of rotations at $v$ which are fixed by $\gamma$ under conjugation. If an
automorphism $\gamma$ fixes a rotation system $\rho$, then $\gamma \rho_v \gamma^{-1}=\rho_{\gamma(v)}$ for each vertex
$v$. It follows that $\gamma^{l(v,\gamma)}\rho_v \gamma^{-l(v,\gamma)}=\rho_v$, or $\rho_v$ is a member of the set
$F_v(\gamma^{l(v,\gamma)})$.  It is the cardinalities of these sets which determine the cardinality of $F(\gamma )$.

\begin{theorem}
\label{thm:fix} If $\gamma \in \Gamma$, then $\abs{F(\gamma)}=\prod_{v \in S}
\abs{F_v\left(\gamma^{l(v,\gamma)}\right)}$, where the product extends over a complete set $S$ of orbit representatives
of $\gamma$ acting on the vertex set $V$.
\end{theorem}

A permutation is defined to be \emph{d-regular} if each of its orbits has cardinality $d$. That is, in the disjoint
cycle representation of the permutation, all the cycles have length $d$. The next result provides the number of
rotations at a vertex which are fixed by the automorphism $\gamma^{l(v,\gamma)}$. We let $\phi$ denote the Euler
function and recall that degree$(v)$ denotes the cardinality of the neighborhood $N(v)$ of $v$.
\begin{theorem}
\label{thm:counting}
\[ \abs{F_v\left(\gamma^{l(v,\gamma)}\right)}=\left\{
  \begin{array}{ll}
    \phi(d)\left(\frac{\mathrm{degree}(v)}{d}-1\right)!d^{\frac{\mathrm{degree}
   (v)}{d}-1} & \mbox{if $\gamma^{l(v,\gamma)}$ is $d$-regular on $N(v)$,} \\
    0 & \mbox{otherwise.}
  \end{array} \right.
 \]
\end{theorem}

We claim that there is much more information contained in the above counting theorems than previously reported.  To
obtain the additional information we supplement Theorem~\ref{thm:Burnside} with a related but more detailed result.

Each automorphism acting on the vertices of the graph carries with it a \emph{cycle type} which records the number of
orbits of a particular length.  For example, $\gamma =(v_1v_2)(v_3v_4)(v_5v_6v_7)$ has two 2-cycles and a 3-cycle.  The
cycle type of any permutation of an $n$-set is encoded as a monomial in the indeterminates $s_1,s_2,\ldots,s_n$, where
the exponent $j_k$ of $s_k$ in the monomial is the number of cycles of length $k$.  The cycle type of the above
permutation is encoded as $s_2^2s_3$.  For simplicity we let $s$ denote the function which assigns to a permutation the
monomial that encodes its cycle type;  that is, $s(\gamma)=\prod s_k^{j_k}$.

The previous theorems show that the number of maps left fixed by a graph automorphism depends considerably on the cycle
type of the automorphism. We are led to define the \emph{imbedding sum} of a graph $G$, denoted $Z(G)$, as the
polynomial whose terms are the cycle-type monomials. Each automorphism $\gamma$ whose cycle type corresponds to the
monomial $\prod s_k^{j_k}$ contributes $\abs{F(\gamma)}$ to the coefficient of this term in the polynomial. Formally,
we have the following.

\begin{equation}
\label{def:Imbedding Sum}  Z(G)=\frac{1}{\abs{\Gamma(G)}}\sum_{\gamma \in \Gamma(G)} \abs{F(\gamma)}s(\gamma).
\end{equation}

We remark that this polynomial is no more difficult to determine using the counting theorems than the number of
unlabeled imbeddings of the given graph.  To illustrate, the number of unlabeled imbeddings of the complete graph $K_n$
as reported in \cite{MRW88} is
\[ \sum_{d|n}\frac{{(n-2)!}^{n/d}}{d^{n/d}(n/d)!}+\sum_{d|(n-1) \atop d \ne 1}
\frac{\phi (d){(n-2)!}^{(n-1)/d}}{n-1}. \] The imbedding sum is found to be
\begin{equation}\label{eqn:Zkn}
Z(K_n)=\sum_{d|n}\frac{{(n-2)!}^{n/d}}{d^{n/d}(n/d)!}s_d^{n/d}+\sum_{d|(n-1) \atop d \ne 1} \frac{\phi
(d){(n-2)!}^{(n-1)/d}}{n-1}s_1s_d^{(n-1)/d}.
\end{equation}

The next theorem together with the counting theorems already presented are the main tools for the applications which
follow.  It is an immediate application of Redfield's lemma \cite{R27}.  The notation used is that of
P\'{o}lya~\cite{P37}.  If $\Gamma$ is a permutation group, then the \emph{cycle index} of $\Gamma$ is the polynomial
\begin{equation}
 Z(\Gamma)=\frac{1}{\abs{\Gamma}}\sum_{\gamma \in \Gamma} s(\gamma).
\end{equation}

The notational similarity between the imbedding sum of a graph and the cycle index of a permutation group is justified
with the following theorem, the main theorem of this article.

\begin{theorem}[Decomposition Theorem]
\label{thm:Decomposition} $Z(G)=\sum Z(\Gamma_M(G))$, where the sum is over the set of unlabeled imbeddings of the
graph $G$.
\end{theorem}

\begin{proof}
Redfield's lemma asserts that if the group $\Gamma$ acts on a set (the set of maps) and $s$ is a function from $\Gamma$
to a ring containing the rationals which is constant on conjugacy classes, then
\[ \frac{1}{\abs{\Gamma}}\sum_{\gamma \in \Gamma} \abs{F(\gamma)}s(\gamma)= \sum \frac{1}{\abs{\Gamma_M}}\sum_{\gamma \in \Gamma_M} s(\gamma),\]
where the sum on the righthand side is over a complete set of orbit representatives (the unlabeled imbeddings). Taking
the function $s$ as the cycle-type monomial we observe that the lefthand side is $Z(G)$ and the righthand side is $\sum
Z(\Gamma_M(G))$.
\end{proof}


\section{P\'{o}lya Enumeration}
\label{section3} The most immediate application of the imbedding sum is its evaluation at $s_k=1$ for all $k$ which is
denoted as $Z(G;1)$. In this case, we have the number of unlabeled imbeddings of the graph. For example, $Z(K_4;1)=3$,
giving three unlabeled imbeddings of $K_4$. One of these is on the sphere and the other two are on the torus.

Suppose instead we have a function which assigns to each vertex of the graph an element of the set $Y$ (say, for
example, a set of colors). Moreover, suppose each element of $Y$ is provided with an $n$-tuple of nonnegative integer
weights. Each assignment is called a \emph{figure}. The \emph{figure counting series} $f(w_1,w_2,\ldots,w_n)$ records
as the coefficient of $w_1^{p_1}w_2^{p_2}\cdots w_n^{p_n}$ the number of members of the set $Y$ that have the $n$-tuple
of weights $(p_1,p_2,\ldots,p_n)$.

Given a particular labeled imbedding $M$ of the graph $G$ and an assignment of members of $Y$ to the vertices of $G$ ,
we describe the assignments by the number of vertices whose assigned members from the set $Y$ have weight
$(p_1,p_2,\ldots,p_n)$.  The map-automorphism group $\Gamma_M$ then acts on these assignments in a natural way to yield
equivalence classes of assignments.  The classes are called \emph{configurations} and are enumerated by weight in the
\emph{configuration counting series}.

The P\'{o}lya Enumeration Theorem~\cite{P37} describes how this series can be determined from the cycle index
$Z(\Gamma_M)$ and the figure counting series $f$. The substitution $s_k^{j_k} \rightarrow
(f(w_1^k,w_2^k,\ldots,w_n^k))^{j_k}$ in the cycle index is denoted $Z(\Gamma_M;f(w_1,w_2,\ldots,w_n))$ and was shown by
P\'{o}lya to provide the configuration counting series (see Harary and Palmer~\cite{HP73} for an excellent exposition
of the theory). In summary, we having the following application of P\'{o}lya's theory.

\begin{theorem}
\label{thm:Map Configurations} If $f(w_1,w_2,\ldots,w_n)$ is the figure counting series for $Y$, then the series which
enumerates by weight the number of configurations of the map $M$ is $Z(\Gamma_M;f(w_1,w_2,\ldots,w_n))$.
\end{theorem}

If we sum these series over a complete set of congruence class representative maps (the unlabeled imbeddings), then we
obtain a series which enumerates by weight the quantity of configurations among the unlabeled imbeddings of the graphs.
Moreover, from the Decomposition Theorem~(\ref{thm:Decomposition}), this sum is the imbedding sum evaluated in the
above manner.  This evaluation is denoted as $Z(G;f(w_1,w_2,\ldots,w_n))$.  We present the following theorem.

\begin{theorem}
\label{thm:Imbedding Configurations} If $f(w_1,w_2,\ldots,w_n)$ is the figure counting series for $Y$, then the series
which enumerates by weight the number of configurations of the graph $G$ among its unlabeled imbeddings is
$Z(G;f(w_1,w_2,\ldots,w_n))$.
\end{theorem}

As an example, if we assign to each vertex of $G$ a color from the set $Y=\{\mathrm{black},\mathrm{white}\}$ whose
members have weight (1,0) and (0,1), respectively, then the figure counting series is $f(b,w)=b^1w^0+b^0w^1=b+w$.  In
the case that $G$ is the complete graph on four vertices we have from Equation~\ref{eqn:Zkn}
\[ Z(K_4)=\frac{1}{24}(16s_1^4+32s_1s_3+12s_2^2+12s_4) \qquad \mbox{and} \]
\begin{eqnarray*}
Z(K_4;b+w)&=&\frac{1}{24}\left(16(b+w)^4+32(b+w)(b^3+w^3)+12(b^2+w^2)^2+12(b^4+w^4)\right)\\
&=&3b^4+4b^3w+5b^2w^2+4bw^3+3w^4.
\end{eqnarray*}

\begin{figure}[h]
\begin{center}
\includegraphics{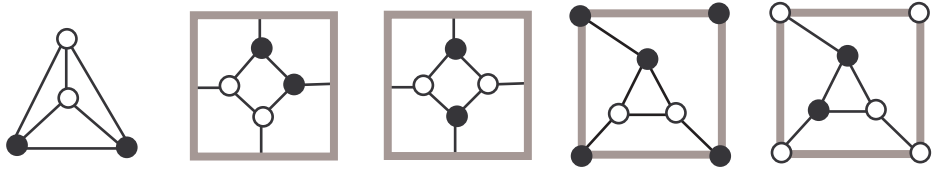}
\caption{The different ways to color two vertices black among the unlabeled imbeddings of $K_4$. The last four
imbeddings are on the torus.}
\label{fig1}
\end{center}
\end{figure}

Evidently there are five different ways to color the vertices in the three unlabeled imbeddings of $K_4$ so that two
vertices receive the color black and two vertices receive the color white (the coefficient 5 of the term $b^2w^2$).
These five colored unlabeled imbeddings are shown in Figure~\ref{fig1}.

Observing that the two black vertices determine a unique edge in $K_4$ we have counted the number of different ways to
root its unlabeled imbeddings at an edge. This technique can be used to determine the edge-rooted imbeddings of any
complete graph, if desired.

In how many ways can the unlabeled imbeddings of $K_n$ be rooted at a vertex?  This question is answered by extracting
the coefficient of $b^1w^{n-1}$ in the configuration counting series $Z(K_n;b+w)$.  This series is obtained from
Equation~\ref{eqn:Zkn} as
\begin{eqnarray}\label{eqn:rootKn}
  \lefteqn{Z(K_n;b+w)=\sum_{d|n}\frac{{(n-2)!}^{n/d}}{d^{n/d}(n/d)!}(b^d+w^d)^{n/d}}\nonumber\\
  &&{}+\sum_{d|(n-1) \atop d \ne 1}\frac{\phi (d){(n-2)!}^{(n-1)/d}}{n-1}(b+w)(b^d+w^d)^{(n-1)/d}.
\end{eqnarray}

Extracting the coefficient of $b^1w^{n-1}$ and simplifying we have that the number of ways to root the unlabeled
imbeddings of $K_n$ at a single vertex is
\begin{equation}\label{E3.2}
  \frac{1}{n-1}\sum_{d|n-1} \phi(d) {(n-2)!}^{(n-1)/d}.
\end{equation}

Curiously, this resembles the cycle index of the cyclic group of order $n-1$.  That is, if $C_{n-1}$ is the cyclic
permutation group of order $n-1$ and degree $n-1$, then
\[ Z(C_{n-1})= \frac{1}{n-1}\sum_{d|n-1} \phi(d) s_d^{(n-1)/d},\]
so that
\begin{equation}\label{E3.3}
  Z(C_{n-1};(n-2)!)= \frac{1}{n-1}\sum_{d|n-1} \phi(d) {(n-2)!}^{(n-1)/d}.
\end{equation}

This result will become less surprising after we explore vertex-rooted graphs in detail in Section~\ref{section4}



\section[Imbedding symmetries]{The enumeration of graph imbeddings by their
symmetries} \label{section4}
\subsection{Introduction}
For many graphs the Decomposition Theorem~(\ref{thm:Decomposition}) enables us to determine precisely how many
imbeddings (both labeled and unlabeled) there are that have a particular map-automorphism group. To see this, consider
the complete graph $K_4$ which has 16 labeled imbeddings and 3 unlabeled imbeddings.

One of the unlabeled imbeddings occurs on the sphere with map-automorphism group $A_4$ (the alternating group of order
12 represented as a permutation group of degree 4).  The index of this group in the graph automorphism group is two so
by Theorem~\ref{thm:index} there are two labeled imbeddings in this congruence class.

The two other unlabeled imbeddings occur on the torus.  One of them has as map-automorphism group the cyclic group of
order 4, $C_4$ (represented as a permutation group of degree 4).  The other toroidal imbedding has the cyclic group of
order 3, $E_1\times C_3$ (also represented as a permutation group of degree 4 where $E_1$ is the identity permutation
group of degree 1).  There are 6 labeled imbeddings in the first congruence class and 8 in the other.

The cycle indexes of these groups and the imbedding sum of $K_4$ satisfy the Decomposition
Theorem~(\ref{thm:Decomposition}) as indicated below.
\[
\begin{array}{lcccccr}
    Z(A_4) &=&\frac{1}{12} &(s_1^4 &+8s_1s_3 &+3s_2^2 &)\\
    Z(C_4) &=& \frac{1}{4} &(s_1^4 & & +s_2^2 &+2s_4)\\
    Z(E_1\times C_3) &=& \frac{1}{3} &(s_1^4 &+2s_1s_3 & &)\\
    \multicolumn{7}{c}{\hrulefill}\\
    Z(K_4) &=& \frac{1}{24} &(16s_1^4&+32s_1s_3 &+12s_2^2 &+12s_4)
\end{array}
\]

Thus, $Z(K_4)$ decomposes as the sum of the cycle indexes $Z(A_4)$, $Z(C_4)$, and $Z(E_1 \times C_3)$.  It is not
difficult to show that no other decomposition of $Z(K_4)$ is possible.  That is, if $Z(K_4)=\sum i_k Z(\Gamma_k)$,
where the $i_k$ are nonnegative integers and the $\Gamma_k$ are permutation groups of degree 4 contained in
$\Gamma(K_4)$, then the only solution is that given above.

Once again applying the P\'{o}lya Enumeration Theory of Section~\ref{section3} with the figure counting series $b+w$ we
have
\begin{equation}
  Z(K_4;b+w)=Z(A_4;b+w)+Z(C_4;b+w)+Z(E_1\times C_3;b+w).
\end{equation}
Expanding each term on the righthand side and extracting the coefficients of $b^2w^2$ with the operator $[b^2w^2]$ we
have
\begin{eqnarray*}
  [b^2w^2]Z(K_4;b+w)&=& [b^2w^2]Z(A_4;b+w)+[b^2w^2]Z(C_4;b+w)+[b^2w^2]Z(E_1\times
  C_3;b+w)\\
  &=& 1+2+2.
\end{eqnarray*}
Thus, among the five unlabeled imbeddings of $K_4$ with two black and two white vertices, one has symmetry group $A_4$,
two have $C_4$, and two have $E_1\times C_3$.  Another glance at Figure~\ref{fig1} confirms this result.

It is unfortunate that uniqueness of the decomposition is not the rule. However, for some graphs many permutation
groups can be eliminated from consideration as map-automorphism groups. In some cases then, additional information
about the graph eliminates the ambiguity in the decomposition. We present some theorems which are useful in this
regard. All but the last of these theorems are well known and their proofs may be found in White~\cite{W84} or Biggs
and White~\cite{BW79}, for example.

\begin{theorem}
\label{thm:Adjacent Vertices} If a map-automorphism group fixes two adjacent vertices, then it is the identity
automorphism.
\end{theorem}

Thus, it is not possible for the dihedral group $D_4$ to be a map-automorphism group of $K_4$ since as a permutation
subgroup of $\Gamma(K_4)$ it would contain a permutation of cycle-type $s_1^2s_2$. This implies that the permutation
fixes two vertices which are necessarily adjacent in $K_4$.

The next theorem will be used extensively to decompose the imbedding sum of wheel graphs, bouquets of loops, and
vertex-rooted graphs.

\begin{theorem}
\label{thm:Cyclic} If each member of a map-automorphism group $\Gamma_M$ fixes the same vertex $v$, then $\Gamma_M$ is
a cyclic permutation group. Moreover, if the neighborhood of $v$ has cardinality degree$(v)$, then a generator of
$\Gamma_M$ restricted to the neighborhood of $v$ is $d$-regular for some divisor $d$ of degree$(v)$. This implies that
each member of $\Gamma_M$ is a regular permutation when restricted to the neighborhood of $v$.
\end{theorem}

Another useful result for excluding groups as map-automorphism groups is the following theorem which limits their
order.

\begin{theorem}
\label{thm:Divisor} If a graph has $e$ edges, then the order of a map-automorphism group must be a divisor of $2e$.
Moreover, if the order of some map-automorphism group attains the value $2e$, then the group acts transitively on the
vertices, edges, and regions of the imbedded graph.
\end{theorem}

In most of the applications to follow we will ensure a unique decomposition of the imbedding sum by restricting the
investigation to graphs whose map-automorphism groups are all cyclic (invoking Theorem~\ref{thm:Cyclic}). In these
cases, the number of unlabeled imbeddings with a particular cyclic symmetry can be explicitly found. The cycle indexes
of these cyclic groups have the form $\frac{1}{d}\sum_{k|d}\phi(k)s_k^{n/k}$ (or $s_1s_k^{n/k}$) (see \cite{HP73}).

\begin{theorem}
\label{thm:Mobius} If the map-automorphism groups of a graph are cyclic, then its imbedding sum uniquely decomposes as
the sum of cyclic cycle indexes. Let $G_n$ be such a graph, the size of which is measured by the integer parameter $n$,
and let $i_d$ be the number of unlabeled imbeddings whose map-automorphism group is cyclic of order $d$. Suppose
further that there is a function, $f$, of $n$ and $d$ alone such that the imbedding sum of $G_n$ satisfies
\[ Z(G_n)=\sum_{d|n} f(n,d)s_d^{n/d} \quad \mbox{(or $s_1s_d^{n/d}$), then}
\]
\begin{equation}\label{eqn:Mobius}
  i_d=d\sum_{k \atop d|k|n} \frac{f(n,k)}{\phi(k)}\mu(k/d) ,
  \quad \mbox{where $\mu$ is the M\"{o}bius function.}
\end{equation}
\end{theorem}
\begin{proof}
From the Decomposition Theorem~(\ref{thm:Decomposition}) we have
\[ \sum_{d|n}f(n,d)s_d^{n/d}=\sum_{k|n}i_k
\frac{1}{k}\sum_{d|k}\phi(d)s_d^{n/d} \qquad \mbox{(or $s_1s_d^{n/d}$).}
\]
There is a unique solution to this system of equations in the unknowns $i_k$ obtained by comparing coefficients of
$s_n$ (or $s_1s_n$) and then, in turn, the coefficients of $s_d^{n/d}$ (or $s_1s_d^{n/d}$) for each successive divisor
$d$ of $n$ beginning with the largest.  Uniqueness of the solution was actually shown to be the case in Redfield's
original article~\cite{R27}.  Extracting the coefficient of $s_d^{n/d}$ (or $s_1s_d^{n/d}$) yields
\[ f(n,d)=\sum_{k \atop d|k|n} i_k \frac{\phi(d)}{k}\quad \mbox{or}
\quad \frac{f(n,d)}{\phi(d)}=\sum_{k \atop d|k|n}\frac{i_k}{k} .
\]
The latter system of equations (one equation for each divisor $d$ of $n$) can be solved explicitly for $i_d/d$ by
M\"{o}bius inversion.  The details of this particular inversion problem can be found in Hall~\cite{H67}.
\end{proof}

\subsection{Wheel graphs.}
We apply the Decomposition Theorem~(\ref{thm:Decomposition}) to a class of graphs where the map-automorphism groups are
known. This is the class of wheel graphs, $W_{n+1}$, described as a cycle of $n$ vertices each joined to a central
vertex. The counting theorems were applied to this class of graphs in \cite{MRW88} to obtain the number of unlabeled
imbeddings for each of its members. From this information it is easy to determine that the imbedding sum for wheels on
five or more vertices is
\begin{equation}\label{eqn:Zwheel}
  Z(W_{n+1})= \left\{
  \begin{array}{ll}
    \frac{1}{2n}\sum_{d|n}\frac{\phi^2(d)}{d}(2d)^{n/d}(\frac{n}{d}-1)!s_1s_d^{n/d}& \mbox{if $n$ is odd,} \\
     \frac{1}{2n}\sum_{d|n}\frac{\phi^2(d)}{d}(2d)^{n/d}(\frac{n}{d}-1)!s_1s_d^{n/d} +2^{n-3}(\frac{n}{2}-1)!s_1s_2^{n/2} & \mbox{if $n$ is even.}
  \end{array} \right.
\end{equation}

To characterize the map-automorphism groups of the wheel graphs we note that each graph automorphism must fix the
central vertex (provided that $n\ge 4$). Thus, if $M$ is a map of $W_{n+1}$ with map-automorphism group $\Gamma_M$,
then an immediate result of Theorem~\ref{thm:Cyclic} is that $\Gamma_M$ must be $E_1\times C_d[E_{n/d}]$ for some
divisor $d$ of $n$, where $E_1$ is the identity permutation group of degree 1 and $C_d[E_{n/d}]$ is the wreath product
of the cyclic permutation group of order $d$ and degree $d$ about the identity permutation group of degree $n/d$. The
cycle index of this group is $\frac{1}{d}\sum_{k|d}\phi(k)s_1s_k^{n/k}$.

Let $i_d$ be the number of unlabeled imbeddings of $W_{n+1}$ which have map-automorphism group $C_d[E_{n/d}]$. We say
that these imbeddings are \emph{$d$-fold symmetric}. Applying the Decomposition Theorem~(\ref{thm:Decomposition}) we
have that
\begin{equation}
  Z(W_{n+1})=\sum_{d|n}i_d \frac{1}{d}\sum_{k|d}\phi(k)s_1s_k^{n/k}.
\end{equation}

The explicit solution obtained by using Theorem~\ref{thm:Mobius} is given below.
\begin{theorem}
The number $i_d$ of unlabeled imbeddings of the wheel $W_{n+1}$ which are $d$-fold symmetric is given by
\[ i_d =d\sum_{k \atop d|k|n} \frac{f(n,k)}{\phi(k)}\mu(k/d),
 \quad \mbox{where}\]
\[ f(n,k)=\left\{
  \begin{array}{ll}
    \frac{\phi^2(k)}{2nk}(2k)^{n/k}(\frac{n}{k}-1)!& \mbox{if $k \ne 2$,} \\
    \frac{\phi^2(k)}{2nk}(2k)^{n/k}(\frac{n}{k}-1)! +2^{n-3}(\frac{n}{2}-1)! & \mbox{if $k=2$.}
  \end{array} \right.\]
\end{theorem}

The number of unlabeled imbeddings of $W_{n+1}$ which are $d$-fold symmetric is given in Table~\ref{table:uWheel} for
the first few values of $n$. From this information and Theorem~\ref{thm:index} we are able to calculate the number of
labeled imbeddings of $W_{n+1}$ which have $d$-fold symmetry. This information is tabulated in
Table~\ref{table:lWheel}.

\begin{table}
\caption{\label{table:uWheel}The number of unlabeled imbeddings of the wheel $W_{n+1}$ which have $d$-fold symmetry.}
\begin{tabular*}{\textwidth}{@{\extracolsep{\fill}}lrrrrr}
    \hline\quad n &4&5&6&7&8\\
    d \\ \hline\hline
    1 &9&76&617&6582&80399\\
    2 &5&0&42&0&479\\
    3 &0&0&5&0&0\\
    4 &2&0&0&0&4\\
    5 &&4&0&0&0\\
    6 &&&2&0&0\\
    7 &&&&6&0\\
    8 &&&&&4\\ \hline
    Total &16&80&666&6588&80886\\ \hline
\end{tabular*}
\end{table}

\begin{table}
\caption{\label{table:lWheel}The number of labeled imbeddings of the wheel $W_{n+1}$ which have $d$-fold symmetry.}
\begin{tabular*}{\textwidth}{@{\extracolsep{\fill}}lrrrrr}
    \hline\quad n &4&5&6&7&8\\
    d \\ \hline\hline
    1 &72&760&7404&92148&1286384\\
    2 &20&0&252&0&3832\\
    3 &0&0&20&0&0\\
    4 &4&0&0&0&16\\
    5 &&8&0&0&0\\
    6 &&&4&0&0\\
    7 &&&&12&0\\
    8 &&&&&8\\ \hline
    Total &96&768&7680&92160&1290240\\ \hline
\end{tabular*}
\end{table}

\subsection{Bouquets}
Another graph whose map-automorphism groups are cyclic is the bouquet of $n$ loops, $B_n$ (represented as a simple
graph by introducing two vertices upon each loop). For this graph, the graph-automorphism group is the permutation
group $E_1 \times S_n[C_2]$ of order $2^nn!$ and degree $2n+1$, where $S_n$ is the symmetric permutation group.

The imbedding sum is (Rieper~\cite{R87})
\begin{equation}\label{eqn:Zbn}
  Z(B_n)=\frac{1}{2n}\sum_{d|n}\frac{\phi(d)(2n/d)!d^{n/d}}{2^{n/d}(n/d)!}s_1s_d^{2n/d}
  +\frac{1}{2n}\sum_{d|n}\sum_{m=0}^{\lfloor\frac{n-1}{2d}\rfloor}\frac{\phi(2d)(n/d)!d^m}{(\frac{n}{d}-2m)!m!}s_1s_{2d}^{n/d}.
\end{equation}

A consequence of Theorem~\ref{thm:Cyclic} is that each map-automorphism group of the bouquet $B_n, n>1$, is one of the
permutation groups $E_1 \times C_d[E_{2n/d}]$ for some divisor $d$ of $2n$.  A map with this map-automorphism group is
again said to be $d$-fold symmetric.

Performing an analysis like that for the wheel graphs we determine the number of unlabeled and labeled imbeddings of
the bouquet which are $d$-fold symmetric. The results of the calculations are given in the theorem below and are
tabulated in Table~\ref{table:uBouquet} and Table~\ref{table:lBouquet}.

\begin{theorem}
The number $i_d$ of unlabeled imbeddings of the bouquet $B_n$ which are $d$-fold symmetric is given by
\[ i_d =d\sum_{k \atop d|k|n} f(n,k)\mu(k/d), \quad \mbox{where}\]
\[
f(n,k)=\left\{
\begin{array}{ll}
    \frac{(2n/k)!(k/2)^{n/k}}{2n(n/k)!} & \quad \mbox{if $k$ is odd,}\\
    \frac{(2n/k)!(k/2)^{n/k}}{2n(n/k)!}+\frac{1}{2n}\sum\limits_{m=0}^{\lfloor \frac{n-1}{k}\rfloor}
    \frac{(2n/k)!(k/2)^m}{(2n/k -2m)!m!} & \quad \mbox{if $k$ is even.}
\end{array}
\right.
 \]
\end{theorem}

\begin{table}
\caption{\label{table:uBouquet}The number of unlabeled imbeddings of the bouquet $B_n$ which have $d$-fold symmetry.}

\begin{tabular*}{\textwidth}{@{\extracolsep{\fill}}lrrrrrr}
    \hline\quad n &1&2&3&4&5&6\\
    d \\ \hline\hline
    1 &0&0&1&10&86&837\\
    2 &1&1&2&5&16&52\\
    3 &&0&1&0&0&5\\
    4&&1&0&2&0&4\\
    5 &&&0&0&2&0\\
    6 &&&1&0&0&3\\
    7 &&&&0&0&0\\
    8 &&&&1&0&0\\
    9 &&&&&0&0\\
    10 &&&&&1&0\\
    11 &&&&&&0\\
    12 &&&&&&1\\ \hline
    Total &1&2&5&18&102&902\\\hline
\end{tabular*}
\end{table}

\begin{table}
\caption{\label{table:lBouquet}The number of labeled imbeddings of the bouquet $B_n$ which have $d$-fold symmetry.}
\begin{tabular*}{\textwidth}{@{\extracolsep{\fill}}lrrrrrr}
    \hline\quad n &1&2&3&4&5&6\\
    d \\ \hline\hline
    1 &0&0&48&3840&330240&38568960\\
    2 &1&4&48&960&30720&1198080\\
    3 &&0&16&0&0&76800\\
    4 &&2&0&192&0&46080\\
    5 &&&0&0&1536&0\\
    6 &&&8&0&0&23040\\
    7 &&&&0&0&0\\
    8 &&&&48&0&0\\
    9 &&&&&0&0\\
    10 &&&&&384&0\\
    11 &&&&&&0\\
    12 &&&&&&3840\\ \hline
    Total &1&6&120&5040&362880&39916800\\ \hline
\end{tabular*}
\end{table}

\subsection{The directed bouquet}
We illustrate here that the counting and decomposition theorems can be used to enumerate graph imbeddings for graphs
that have additional structure such as a directed graph. We choose the directed bouquet of $n$ loops,
$\overrightarrow{B_n}$, which is the bouquet with each edge (loop) given an orientation (indicated by placing an arrow
along each loop).

The counting theorems listed in Section~\ref{section3} arose from an analysis of the action of the graph-automorphism
group on the \emph{vertices} of the graph but the graph-automorphism group of the directed bouquet acts on its set of
directed edges. We overcame a similar problem for the undirected bouquet by twice subdividing each loop to yield a
simple graph. This same technique is used here except that one of the two vertices introduced upon an edge is colored
black and the other white. We agree then, that the orientation of the edge proceeds from the black to the white vertex.
The graph automorphisms must map a black vertex to a black vertex, a white vertex to another white vertex, and the
original central vertex to itself.

With this model we can use the counting and decomposition theorems as they were presented with the following
modification. We record the cycle type of any automorphism with respect to its action on the directed edges instead of
its action on the vertices.

The graph-automorphism group of $\overrightarrow{B_n}$ is the symmetric group $S_n$.  Theorem~\ref{thm:counting}
implies that we only need to consider the permutations in $S_n$ that are $d$-regular for some divisor $d$ of $n$.  Let
$\gamma$ be one of the $\frac{n!}{d^{n/d}(n/d)!}$ $d$-regular permutations in $S_n$.  Let $v$ be the central vertex,
then $l(v,\gamma)=1$ so the fixed set at $v$ satisfies
\[ \abs{F_v(\gamma)}=\phi(d)\left(\frac{2n}{d}-1\right)!d^{\frac{2n}{d}-1}.\]

If $u$ is a noncentral vertex (either black or white), then it has two neighbors one of which is the central vertex.
Thus, $\gamma^{l(u,\gamma)}$ is 1-regular on the neighborhood of $u$ so by Theorem~\ref{thm:counting} the fixed set at
$u$ has cardinality 1.

From this information and Theorem~\ref{thm:Burnside} the imbedding sum of $\overrightarrow{B_n}$ is given by
\begin{eqnarray}\label{eqn:Bndirected1}
  Z(\overrightarrow{B_n})
  &=&\frac{1}{n!}\sum_{d|n}\frac{n!}{d^{n/d}(n/d)!}
  \phi(d)\left(\frac{2n}{d}-1\right)!d^{\frac{2n}{d}-1}s_d^{n/d}\\
  &=&\sum_{d|n}\frac{\phi(d)(2n/d)!d^{n/d}}{2n(n/d)!}s_d^{n/d}.\nonumber
\end{eqnarray}

The map-automorphism groups are once again forced to be cyclic by the presence of the central vertex. Each must be
$C_d[E_{n/d}]$ for some divisor $d$ of $n$ with cycle index $\frac{1}{d}\sum_{k|d}\phi(k)s_k^{n/k}$. Letting $i_d$
denote the number of unlabeled imbeddings with this map-automorphism group we have from
Theorem~(\ref{thm:Decomposition}) that
\begin{equation}\label{eqn:Bndirected2}
  Z(\overrightarrow{B_n})=\sum_{d|n}i_d\frac{1}{d}\sum_{k|d}\phi(k)s_k^{n/k}.
\end{equation}
Equating the righthand sides of Equation~\ref{eqn:Bndirected1} and Equation~\ref{eqn:Bndirected2} and using
Theorem~\ref{thm:Mobius} we have
\begin{theorem}
The number $i_d$ of unlabeled imbeddings of the directed bouquet $\overrightarrow{B_n}$ which are $d$-fold symmetric is
given by
\begin{equation}\label{eqn:Bndirected3}
  i_d=\frac{d}{2n}\sum_{k \atop
  d|k|n}\mu(k/d)\frac{(2n/k)!k^{n/k}}{(n/k)!}.
\end{equation}
\end{theorem}

Some results are tabulated in Table~\ref{table:uBndirect} and Table~\ref{table:lBndirect}. Note that the number of
unlabeled imbeddings of $\overrightarrow{B_n}$ with the largest map-automorphism group, $C_n$,
 is $n$. This is easily confirmed from Equation~\ref{eqn:Bndirected3} above.

\begin{figure}
\begin{center}
\includegraphics[scale=.75, keepaspectratio=true]{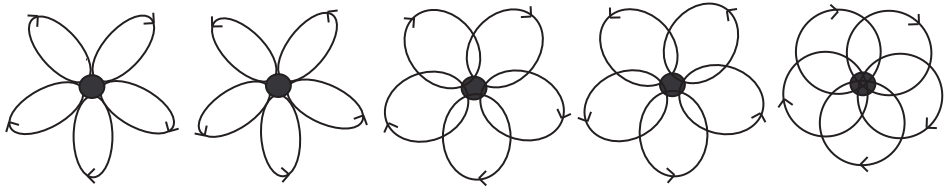}
\caption{The five different imbeddings of the directed bouquet $\vec{B_5}$ of maximum symmetry represented as overlap
graphs.  The first two drawings represent planar imbeddings and the last three represent toroidal imbeddings.}
\label{fig2}
\end{center}
\end{figure}

In Figure~\ref{fig2} we depict the five unlabeled imbeddings of $\overrightarrow{B_5}$ as \emph{overlap} graphs. The
overlap graphs are obtained from the surface imbeddings by projection onto the plane. These overlap graphs are very
useful constructs and in our case capture the 5-fold symmetry of the imbeddings very nicely.  However, these drawings
can be misleading.  Why isn't there a sixth drawing paired with the fifth by a reversal of the directed edges as is the
case for the first two pairs?  We leave the answer to the curious reader.

\begin{table}
\caption{\label{table:uBndirect}The number of unlabeled imbeddings of the directed
 bouquet $\vec{B_n}$ which have $d$-fold symmetry.}
\begin{tabular*}{\textwidth}{@{\extracolsep{\fill}}lrrrrrrrr}
    \hline\quad n &1&2&3&4&5&6&7&8\\
    d \\ \hline\hline
    1 &1&3&20&204&3023&55352&1235519&32430720\\
    2 &&2&0&10&0&158&0&3336\\
    3 &&&3&0&0&24&0&0\\
    4 &&&&4&0&0&0&44\\
    5 &&&&&5&0&0&0\\
    6 &&&&&&6&0&0\\
    7 &&&&&&&7&0\\
    8 &&&&&&&&8\\ \hline
    Total &1&5&23&218&3028&55540&1235526&32434108\\\hline
\end{tabular*}
\end{table}

\begin{table}
\caption{\label{table:lBndirect}The number of labeled imbeddings of the directed bouquet $\vec{B_n}$ which have
$d$-fold symmetry.}
\begin{tabular*}{\textwidth}{@{\extracolsep{\fill}}lrrrrrr}
    \hline\quad n &1&2&3&4&5&6\\
    d \\ \hline\hline
    1 &1&6&120&4896&362760&39853440\\
    2 &&2&0&120&0&56880\\
    3 &&&6&0&0&5760\\
    4 &&&&24&0&0\\
    5 &&&&&120&0\\
    6 &&&&&&720\\ \hline
    Total&1&8&126&5040&362880&39916800\\ \hline
\end{tabular*}
\end{table}

\subsection{Vertex-rooted graphs}
As mentioned previously, if we desire to enumerate the imbeddings of a graph by their map-automorphism groups, then we
need more information beyond its imbedding sum. Information about its map-automorphism groups must be available because
the decomposition of the imbedding sum given by the Decomposition Theorem~(\ref{thm:Decomposition}) is not, in general,
unique. In the examples given so far, the map-automorphism groups were all cyclic and in this case the decomposition of
the imbedding sum can be found.

If the map-automorphism groups of the graph are not all cyclic, then we have two choices.  Either additional
information must be made available or we must alter the graph to force the map-automorphism groups to be cyclic.  The
latter can be accomplished by rooting the graph at a vertex.  In general,  there are numerous choices for the root but
we denote any one of them as $G^*$.  Every automorphism of the rooted graph must fix the root and, hence, by
Theorem~\ref{thm:Cyclic}, each map-automorphism group is cyclic.  Information about the imbedding symmetries of the
rooted graph is then available as in the previous cases.

As an example, we consider the vertex-rooted complete graph $K_n^*$. The automorphism group of this graph is $E_1
\times S_{n-1}$. Since each automorphism must fix the rooted vertex, Theorem~\ref{thm:Adjacent Vertices} implies that
no non-identity map-automorphism can fix another vertex and Theorem~\ref{thm:Cyclic} implies we need only consider
those automorphisms of cycle-type $s_1s_d^{(n-1)/d}$ for some divisor $d$ of $n-1$. There are $(n-1)!/(d^{(n-1)/d}\cdot
((n-1)/d)!)$ such permutations in $E_1\times S_{n-1}$. Let $\gamma$ be one of these permutations and let $u$ be a
vertex other than the root. Then $l(u,\gamma)=d$ so that $\gamma^{l(u,\gamma)}$ is the identity permutation. Applying
Theorem~\ref{thm:counting}, the fixed set at $u$ for $\gamma$ satisfies $ \abs{F_u(\gamma^{l(u,\gamma)})}=(n-2)!$ since
$\gamma^{l(u,\gamma)}$ is 1-regular (it is the identity). There are $\frac{n-1}{d}$ orbit representatives like the
vertex $u$ under the action of $\gamma$ so their contribution to the cardinality of the fixed set of $\gamma$ is the
factor $(n-2)!^{(n-1)/d}$ (see Theorem~\ref{thm:fix}).

If $v$ is the rooted vertex, then $l(v,\gamma)=1$ and we have from Theorem~\ref{thm:counting} that
\[ \abs{F_v(\gamma^{l(v,\gamma)})}=\phi(d)\left(\frac{n-1}{d}-1\right)!d^{\frac{n-1}{d}-1}.\]
Thus, the fixed set for $\gamma$ satisfies
\[ \abs{F(\gamma)}=\phi(d)\left(\frac{n-1}{d}-1\right)!d^{\frac{n-1}{d}-1}(n-2)!^{\frac{n-1}{d}}.\]
Summing over all the $d$-regular permutations we find after simplifying that the imbedding sum is
\begin{equation}\label{eqn:ZK*n}
  Z(K_n^*)=\frac{1}{n-1}\sum_{d|(n-1)}\phi(d){(n-2)!}^{(n-1)/d}s_1s_d^{(n-1)/d}.
\end{equation}
This is in agreement with Equation~\ref{E3.2} derived by applying the P\'{o}lya Enumeration Theory to the imbedding sum
of the unrooted complete graph.

Each of the map-automorphism groups of $K_n^*$ must be $E_1\times C_d[E_{(n-1)/d}]$ for some divisor $d$ of $n-1$, in
which case the corresponding map is $d$-fold symmetric. The imbedding sum uniquely decomposes as a sum of these groups.
Applying Theorem~\ref{thm:Mobius} we have
\begin{theorem}
The number $i_d$ of unlabeled imbeddings of the vertex-rooted complete graph $K_n^*$ which are $d$-fold symmetric is
given by
\[ i_d =d\sum_{k \atop d|k|n-1}\mu(k/d)\frac{(n-2)!^{(n-1)/k}}{n-1}.\]
\end{theorem}

The numbers of unlabeled and labeled imbeddings of $K_n^*$ which are $d$-fold symmetric are given in
Table~\ref{table:uK*n} and Table~\ref{table:lK*n}, respectively.

\begin{table}
\caption{\label{table:uK*n}The number of unlabeled imbeddings of the vertex-rooted complete graph $K_n^*$ which have
$d$-fold symmetry.}
\begin{tabular*}{\textwidth}{@{\extracolsep{\fill}}lrrrrrrr}
    \hline\quad n &1&2&3&4&5&6&7\\
    d \\ \hline\hline
    1 &1&1&0&2&315&1592520&497662709620\\
    2 &&0&1&0&15&0&575960\\
    3 &&&0&2&0&0&7140\\
    4 &&&&0&6&0&0\\
    5 &&&&&0&24&0\\
    6 &&&&&&0&120\\ \hline
    Total &1&1&1&4&336&1592548&497663292840\\\hline
\end{tabular*}
\end{table}

\begin{table}
\caption{\label{table:lK*n}The number of labeled imbeddings of the vertex-rooted complete graph $K_n^*$ which have
$d$-fold symmetry.}
\begin{tabular*}{\textwidth}{@{\extracolsep{\fill}}lrrrrrr}
    \hline\quad n &1&2&3&4&5&6\\
    d \\ \hline\hline
    1 &1&1&0&12&7560&191102400\\
    2 &&0&1&0&180&0\\
    3 &&&0&4&0&0\\
    4 &&&&0&36&0\\
    5 &&&&&0&576\\
    6 &&&&&&0\\ \hline
    Total &1&1&1&16&7776&191102976\\ \hline
\end{tabular*}
\end{table}

Of special interest is the number of imbeddings which have the largest possible map-automorphism group, $E_1\times
C_{n-1}$.  The number of unlabeled imbeddings is $i_n=(n-2)!$.  There are $(n-2)!$ different labeled imbeddings in each
of these congruence classes.  In the case $n$ equals 4, two labeled imbeddings on the sphere and two on the torus have
3-fold symmetry.

\subsection{The complete graph $K_5$.}
The alternative to rooting the graph at a vertex is additional information about its map-automorphism groups. We
illustrate by decomposing the imbedding sum of the complete graph $K_5$. There are $(4!)^5=7776$ different labeled
imbeddings and $Z(K_5;1)=78$ unlabeled imbeddings.  The imbedding sum (Equation~\ref{eqn:Zkn}) is
\begin{equation}
  Z(K_5)=\frac{1}{120}(7776s_1^5+1080s_1s_2^2+360s_1s_4+144s_5).
\end{equation}

There are numerous decompositions of this sum as the sum of cycle indexes of permutation groups.  However, from
Theorem~\ref{thm:Divisor} we know that each map-automorphism group of $K_5$ must have an order that is a divisor of 20.
Thus, for example, we can exclude the groups $E_1\times C_2[C_2]$, $E_2\times C_3$, and $C_2\times C_3$.
Theorem~\ref{thm:Adjacent Vertices} implies that we can also exclude $E_1\times C_2\times C_2$ and others which have a
nonidentity element that fixes two adjacent vertices.

In addition, it is known that a complete graph with a prime power number of vertices has a map-automorphism group of
maximum order and that each is a Frobenius group (Biggs~\cite{B71}).  For $K_5$ we denote this group as $F$ and report
that its cycle index is $Z(F)=\frac{1}{20}(s_1^5+5s_1s_2^2+10s_1s_4+4s_5)$.

The other possible map-automorphism groups are the dihedral group $D_5$, the identity group $E_5$, and the cyclic
groups $C_5$, $E_1\times C_4$, and $E_1\times C_2[E_2]$. We let $i_F$ be the number of unlabeled imbeddings with
map-automorphism group the Frobenius group, $i_5$, $i_4$, and $i_2$ the quantity with cyclic map-automorphism groups
$C_5$, $E_1 \times C_4$, and $E_1 \times C_2[E_2]$, respectively, $i_D$ the number with dihedral symmetry, and $i_1$
the number of asymmetric unlabeled imbeddings (symmetry group $E_5$). We are led to solve the equation

\begin{eqnarray}
\label{eqn:ZK5} \lefteqn{\frac{1}{120}(7776s_1^5+1080s_1s_2^2+360s_1s_4+144s_5)=}\hspace{1in}\\ & &
\begin{array}{cl@{\cdot}rrcccl}
   &i_F &\frac{1}{20}(&s_1^5 &+5s_1s_2^2 &+10s_1s_4 &+4s_5&)\\
  +&i_5 &\frac{1}{5} (&s_1^5 & & &+4s_5&)\\
  +&i_4 &\frac{1}{4} (&s_1^5 &+s_1s_2^2 &+2s_1s_4 & &)\\
  +&i_2 &\frac{1}{2} (&s_1^5 &+s_1s_2^2 & & &)\\
  +&i_D &\frac{1}{10}(&s_1^5 &+5s_1s_2^2 & &+4s_5&)\\
  +&i_1 &\frac{1}{1} (&s_1^5 & & &&).
\end{array} \nonumber
\end{eqnarray}

By comparing coefficients there are four different solutions to the above.  One of the solutions gives $i_5$ equal to
one, the other three solutions give $i_5$ equal to zero.  Thus, if we can construct a map of $K_5$ which has $C_5$ as
its map-automorphism group, then we would be done.  This we do now.

Label the vertices of $K_5$ with the integers 1 through 5 and let $\gamma$ be the automorphism $(1,2,3,4,5)$.  Choose
the rotation at vertex 1 to be the cyclic permutation $(2,3,4,5)$ and let $\gamma$ act on this permutation repeatedly
by conjugation to produce the rotations at the other vertices.  The rotation system is then $\rho_1=(2,3,4,5)$,
$\rho_2=(3,4,5,1)$, $\rho_3=(4,5,1,2)$, $\rho_4=(5,1,2,3)$, and $\rho_5=(1,2,3,4)$.

The map-automorphism group of the corresponding imbedding then has $C_5$ as a subgroup.  From the rotation system we
find that the imbedding has two 5-sided regions and one 10-sided region (see Gross and Tucker~\cite{GT87} for the
algorithm which produces the region sizes from the rotation system).  Thus, we know from Theorem~\ref{thm:Divisor} that
the map-automorphism group can not be the Frobenius group $F$ which requires all regions to be of the same size.  The
group is then either $C_5$ or the dihedral group $D_5$.  If it is the latter, then the automorphism
$\alpha=(1)(2,5)(3,4)$ must fix the rotation system above.  However $\alpha$ takes $\rho_1=(2,3,4,5)$ to $(5,4,3,2)$ by
conjugation and this is not $\rho_{\alpha 1}=\rho_1$.

Thus, we have produced an imbedding of $K_5$ which has map-automorphism group $C_5$. The decomposition of the imbedding
sum is now uniquely determined to be $i_F=2$, $i_5=1$, $i_D=0$, $i_4=4$, $i_2=15$ and $i_1=56$.

\end{document}